\documentclass[preprint,12pt]{elsarticle}

\usepackage{amssymb}
\usepackage{amsmath}
\usepackage{amsthm}
\newtheorem{definition}{Definition}

\journal{J. Comput. Appl. Math}

\begin{document}

\begin{frontmatter}
\title{Orthogonal polynomials for Minkowski's question mark function}
\author{Zo\'e Dresse}
\author{Walter Van Assche\corref{cor}}
\cortext[cor]{Corresponding author}
\ead{walter@wis.kuleuven.be}
\address{Department of Mathematics, KU Leuven, Celestijnenlaan 200B box 2400, BE-3001 Leuven, Belgium}

\begin{abstract}
Hermann Minkowski introduced a function in 1904 which maps quadratic irrational numbers to rational numbers and this function is now known as Minkowski's question mark function since Minkowski used the notation $?(x)$. 
This function is a distribution function on $[0,1]$ which defines a singular continuous measure with support $[0,1]$. 
Our interest is in the (monic) orthogonal polynomials $(P_n)_{n \in \mathbb{N}}$ for the Minkowski measure and in particular in the behavior of the recurrence coefficients of the three term recurrence relation.  We will give some numerical experiments using the discretized Stieltjes-Gautschi
method with a discrete measure supported on the Minkowski sequence. We also explain how one can compute the moments
of the Minkowski measure and compute the recurrence coefficients using the Chebyshev algorithm. 
\end{abstract}

\end{frontmatter}


\section{Introduction}   \label{sec:intro}
In 1904 Hermann Minkowski \cite{Minkowski} introduced an interesting function, which he called the question mark function and he denoted
its values by $?(x)$. This notation with a question mark is somewhat confusing, so instead we will denote the function by $q$ and we will 
only consider it on the interval $[0,1]$. 

There are several ways to define the Minkowski question mark function. Minkowski used the following construction: let $\mathcal{M}_1$ be the
sequence with two elements $0$ and $1$ and define $q(0)=0$ and $q(1)=1$. The sequence $\mathcal{M}_2$ then consists of $\mathcal{M}_1$ and the new point
$(0+1)/(1+1)=1/2$ and $q(1/2)=1/2$. In general we construct the \textit{Minkowski sequence} $\mathcal{M}_N$ by taking all the elements from $\mathcal{M}_{N-1}$ and all the
\textit{``mediants''} $(a+a')/(b+b')$ of two consecutive rational numbers $a/b$ and $a'/b'$ in $\mathcal{M}_{N-1}$, where we take $0=0/1$ and $1=1/1$.
Then the Minkowski question mark function on the new points takes the values
\[    q(\frac{a+a'}{b+b'}) = \frac{q(a/b) + q(a'/b')}{2}.  \]
The Minkowski sequence $\mathcal{M}_N$ is dense in $[0,1]$ as $N \to \infty$ and $q(x)$ for $x \in [0,1] \setminus \mathbb{Q}$ is defined by continuity.
Observe that $\mathcal{M}_N$ contains $2^{N-1}+1$ points.

Another way to define the question mark function is by using continued fractions \cite{Denjoy1}. 
If $0 < x < 1$ then we can write $x$ as a regular continued fraction
\[    x = \cfrac{1}{a_1+\cfrac{1}{a_2+\cfrac{1}{a_3 + \cfrac{1}{\ddots}}}}, \qquad a_i \in \mathbb{N} \setminus \{0\}. \]
The Minkowski question mark function at $x$ is then defined as
\[    q(x) = 2 \sum_{k=1}^\infty \frac{(-1)^{k+1}}{2^{a_1+a_2+\cdots +a_k}}. \]
If $x$ is a rational number, then the continued fraction is terminating and $q(x)$ is given by a finite sum.
By setting $q(0)=0$ and $q(1)=1$ one can show that $q: [0,1] \to [0,1]$ is a continuous and increasing function, so that
$q$ is a probability distribution function on $[0,1]$ which defines a probability measure on $[0,1]$. Arnaud Denjoy \cite{Denjoy2} showed that this distribution function has the property that $q'(x)=0$ almost everywhere on $[0,1]$ so that the corresponding measure is singular and continuous.

A third way is to define the question mark function as a fixed point of an iterated function system with two rational functions. One has
\begin{equation}  \label{self-similar}
     q(x) = \begin{cases}  \frac12 q\left( \frac{x}{1-x} \right), & 0 \leq x \leq \frac12, \\
                                  1 - \frac12 q \left( \frac{1-x}{x} \right), & \frac12 < x \leq 1,
                              \end{cases}  
\end{equation}
and one can easily show that the sequence of probability distribution functions $(q_n)_{n \in \mathbb{N}}$, with 
\[    q_n(x) = \begin{cases}  \frac12 q_{n-1}\left( \frac{x}{1-x} \right), & 0 \leq x \leq \frac12, \\
                                  1 - \frac12 q_{n-1} \left( \frac{1-x}{x} \right), & \frac12 < x \leq 1,
                              \end{cases}  \]
and $q_0$ any probability distribution on $[0,1]$, converges uniformly to Min\-kowski's question mark function. This allows us to compute integrals by a limit procedure
\[   \int_0^1 f(x) \, dq(x) = \lim_{n \to \infty} \int_0^1 f(x) \, dq_n(x).  \]
 
\begin{figure}[ht]
\centering
\rotatebox{270}{\resizebox{3in}{!}{\includegraphics{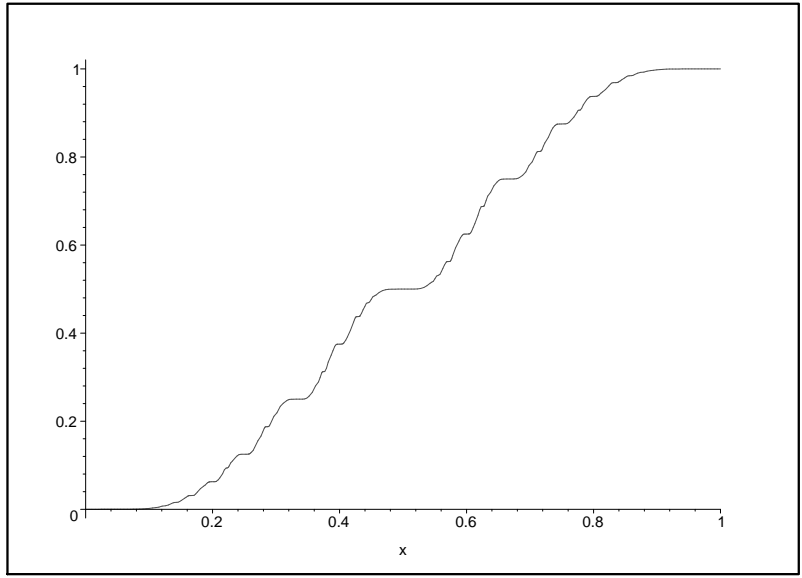}}}
\caption{The Minkowski question mark function}
\label{figq}
\end{figure}

In 1943 Rapha\"el Salem posed a problem about the Fourier coefficients of Minkowski's question mark function:
\[    \alpha_n = \int_0^1 e^{2in\pi x} \, dq(x) . \]
The Riemann-Lebesgue lemma  tells us that Fourier coefficients of an absolutely continuous measure on $[0,1]$ tend to zero.
The Minkowski question mark function is singularly continuous, so one cannot use the Riemann-Lebesgue lemma. Nevertheless, the
support of $q$ is the full interval $[0,1]$ and it was proved by Salem \cite{Salem} that $q$ is H\"older continuous of order
$\alpha=\log 2/(2\log \frac{\sqrt{5}+1}{2})=0.7202$. Furthermore, Salem showed that
\[   \lim_{n \to \infty} \frac{1}{n} \sum_{k=0}^n |\alpha_k| = \mathcal{O}(n^{-\alpha/2}), \]
so that $\alpha_n$ converges to zero on the average and there is the possibility that $\alpha_n \to 0$. 
This is the problem posed by Rapha\"el Salem \cite{Salem}: do the Fourier coefficients of the Minkowski
question mark function converge to 0? 
This is still an open problem. Giedrius Alkauskas  \cite{Alkauskas0} \cite{Alkauskas0.5} already investigated this extensively by both numerical 
and analytical methods.

Our interest in this paper is in the orthonormal polynomials for the Minkowski question mark function:
\[  \int_0^1 p_n(x)p_m(x)\, dq(x) = \delta_{m,n}, \]
where $p_n(x) = \gamma_n x^n + \cdots$ and $\gamma_n > 0$, 
with recurrence relation
\begin{equation}   \label{RR}
  xp_n(x) = a_{n+1} p_{n+1}(x) + b_n p_n(x) + a_n p_{n-1}(x), \qquad n \geq 0, 
\end{equation}
with $p_0=1$ and $p_{-1}=0$, and in particular we are interested in the asymptotic behavior of the recurrence coefficients $(a_n)_{n \geq 1}$ and $(b_n)_{n \geq 0}$.
Rakhmanov's theorem \cite{Rakhmanov} \cite{Rakhmanov2} tells us that for an absolutely continuous measure $\mu$ on $[0,1]$ for which $\mu'>0$ almost everywhere on $[0,1]$, one has
$a_n \to 1/4$ and $b_n \to 1/2$ as $n \to \infty$. In our case $q' = 0$ almost everywhere, so one cannot use Rakhmanov's theorem to deduce
the asymptotic behavior of the recurrence coefficients. However, it is known (see, e.g., \cite{Lubinsky,Magnus,Totik})
that there exist discrete measures and continuous singular measures on $[0,1]$ for which the recurrence coefficients have the behavior $b_n \to 1/2$ and $a_n \to 1/4$ as $n \to \infty$, so that they are in the Nevai class $M(\frac12,\frac14)$. 
\begin{definition}
The Nevai class $M(b,a)$ consists of all positive measures on the real line for which the orthogonal polynomials have recurrence coefficients
satisfying
\[   \lim_{n \to \infty} a_n = a, \qquad \lim_{n \to \infty}  b_n = b.   \]
\end{definition}
It is well known that measures $\mu \in M(b,a)$ have essential spectrum $[b-2a,b+2a]$, i.e., the support of $\mu$ is $[b-2a,b+2a] \cup E$, where
$E$ is at most countable and the accumulation points can only be at $b\pm 2a$ (Blumenthal's theorem, see, e.g., 
\cite[Thm. 7 on p. 23]{Nevai}, \cite[\S 5]{WVA}).

Our first problem is to find out whether
the Minkowski question mark function is such a singular continuous function for which the recurrence coefficients are in the Nevai class (for the interval $[0,1]$), i.e., is the following asymptotic behavior true
\[ \lim_{n \to \infty} a_n = \frac14,  \quad \lim_{n \to \infty} b_n = \frac12 \ ? \] 
The symmetry of $q$ around the point $1/2$
\[    q(x) = 1- q(1-x), \qquad x \in [0,1], \]
already implies that $b_n=1/2$ for all $n \in \mathbb{N}$, so the main problem is to find the asymptotic behavior of the recurrence coefficients $(a_n)_{n \geq 1}$. We will investigate this numerically. In Section 2 we will use the discretized
Stieltjes-Gautschi method to compute the recurrence coefficients by approximating the question mark function by a discrete measure which is the
empirical distribution function of the Minkowski sequence $\mathcal{M}_N$ for various values of $N$. In Section 3 we will use the moments of the
question mark function to compute the recurrence coefficients, using the Chebyshev algorithm. In Section 4 we will compare both methods and our conclusion is that the orthogonal polynomials for the question mark function are probably not in Nevai's class for the interval $[0,1]$.

A second problem is whether the Minkowski question mark function induces a \textit{regular measure} on $[0,1]$ in the sense of Ullman-Stahl-Totik.
Regular measures in the theory of general orthogonal polynomials are those measures for which the asymptotic zero distribution and the $n$th root asymptotics of the leading coefficient $\gamma_n$ of the orthonormal polynomial $p_n$ is given in terms of the equilibrium measure and the capacity of the support $S_\mu$ \cite[Def. 1.7 on p. 123]{Ullman}, \cite[Def. 3.1.2 on p. 61]{StahlTotik}.
\begin{definition}
A positive Borel measure $\mu$ on the real line with compact support $S_\mu$ is a regular measure if
\[   \lim_{n \to \infty}  \gamma_n^{-1/n} = \textup{cap}(S_\mu),  \]
where $\textup{cap}(S_\mu)$ is the logarithmic capacity of $S_\mu$, 
or, equivalently, the zeros $x_{1,n} < x_{2,n} < \cdots < x_{n,n}$ of $p_n$ have the behavior 
\[   \lim_{n \to \infty} \frac1n \sum_{k=1}^n f(x_{k,n})  = \int_{S_\mu} f(x)\, d\mu_e(x)  \]
for every continuous function on $S_\mu$,  where $\mu_e$ is the equilibrium
measure for the set $S_\mu$.
\end{definition}
The capacity of an interval $[a,b]$ is given by $(b-a)/4$, so our problem is: 
does the following asymptotic behavior hold
\[ \lim_{n \to \infty} \gamma_n^{-1/n} = \frac14\ ? \] 
By comparing the coefficient of $x^{n+1}$ in the recurrence relation \eqref{RR} one finds the well-known relation $a_{n+1} = \gamma_n/\gamma_{n+1}$, so that
\[     a_1a_2a_3\cdots a_n = \frac{\gamma_0}{\gamma_n},  \]
and the problem then is to find whether the geometric mean of the recurrence coefficients converges to $1/4$:
\[   \lim_{n \to \infty}  (a_1a_2\cdots a_n)^{1/n} = \frac14 \ ? \]
Of course, when $a_n \to 1/4$, i.e., when
the recurrence coefficients are in the Nevai class for the interval $[0,1]$, then the geometric mean also converges to $1/4$. However, as we mentioned earlier, the numerical experiments in Sections 2 and 3 indicate that the recurrence coefficients probably do not converge to $1/4$, but then it is
still possible that the geometric mean converges to $1/4$. Our numerical results in Section 3 however indicate that this is not the case and that
the geometric mean seems to converge to a value less than $1/4$. 

\section{The discretized Stieltjes-Gautschi method}   \label{Stieltjes-Gautschi}
The computation of the recurrence coefficients of the orthogonal polynomials with the question mark function requires that we need to
be able to integrate polynomials using the measure induced by $q$. This is not easy since the Minkowski question mark function is either
defined by a limiting process or by a series involving continued fraction coefficients. We therefore will compute approximate values of the recurrence coefficients using a discretized version of the Stieltjes method, which was introduced by W. Gautschi, see, e.g., \cite[\S 2.2]{Gautschi}, 
\cite{Gautschi1982}. 
The idea is to use a discrete distribution function $q_N$ (with finitely many points of increase) which converges weakly to the Minkowski question mark function and to compute the recurrence coefficients $(a_{n,N})_{n\geq 1}$ and $(b_{n,N})_{n \geq 0}$ for the orthogonal polynomials for the distribution
$q_N$. Then
\begin{equation}  \label{discreteapprox}
   \lim_{N \to \infty} a_{n,N} = a_n, \quad  \lim_{N \to \infty} b_{n,N} = b_n,  
\end{equation}
where $(a_n)_{n\geq 1}$ and $(b_n)_{n \geq 0}$ are the recurrence coefficients of the orthogonal polynomials for the limiting distribution $q$,
so that the recurrence coefficients of the discrete orthogonal polynomials are approximations of the recurrence coefficients of the orthogonal polynomials for the question mark function. The main advantage is that the computations for the discrete measure only require matrix computations
and can therefore be easily done.

Rather than using the orthonormal polynomials, we will be using the monic orthogonal polynomials $P_n(x) = p_n(x)/\gamma_n$. These
monic orthogonal polynomials satisfy the recurrence relation
\begin{equation}   \label{RRmonic}
    P_{n+1}(x) = (x-b_n) P_n(x) - a_n^2 P_{n-1}(x), \qquad  n \geq 0,  
\end{equation}
with $P_0=1$ and $P_{-1}=0$. In particular we will compute the squared recurrence coefficients $(a_n^2)_{n \geq 1}$.
There is no need to compute the recurrence coefficients $(b_n)_{n \geq 0}$ since the symmetry of $q$ around $1/2$ implies that
$b_n = 1/2$ for all $n \geq 0$. 

\begin{figure}[h!t]
\centering
\rotatebox{270}{\resizebox{2.3in}{!}{\includegraphics{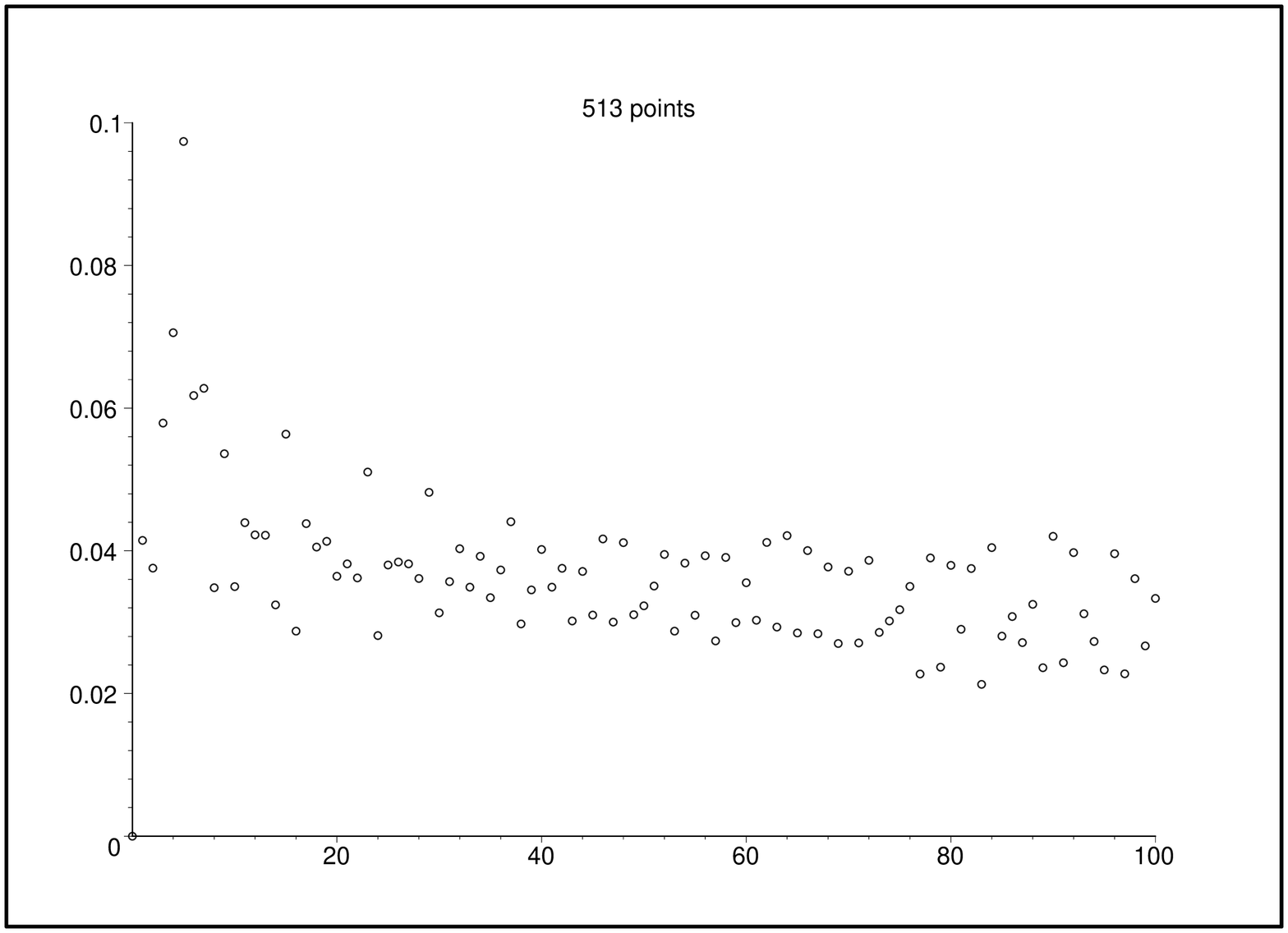}}}
\rotatebox{270}{\resizebox{2.3in}{!}{\includegraphics{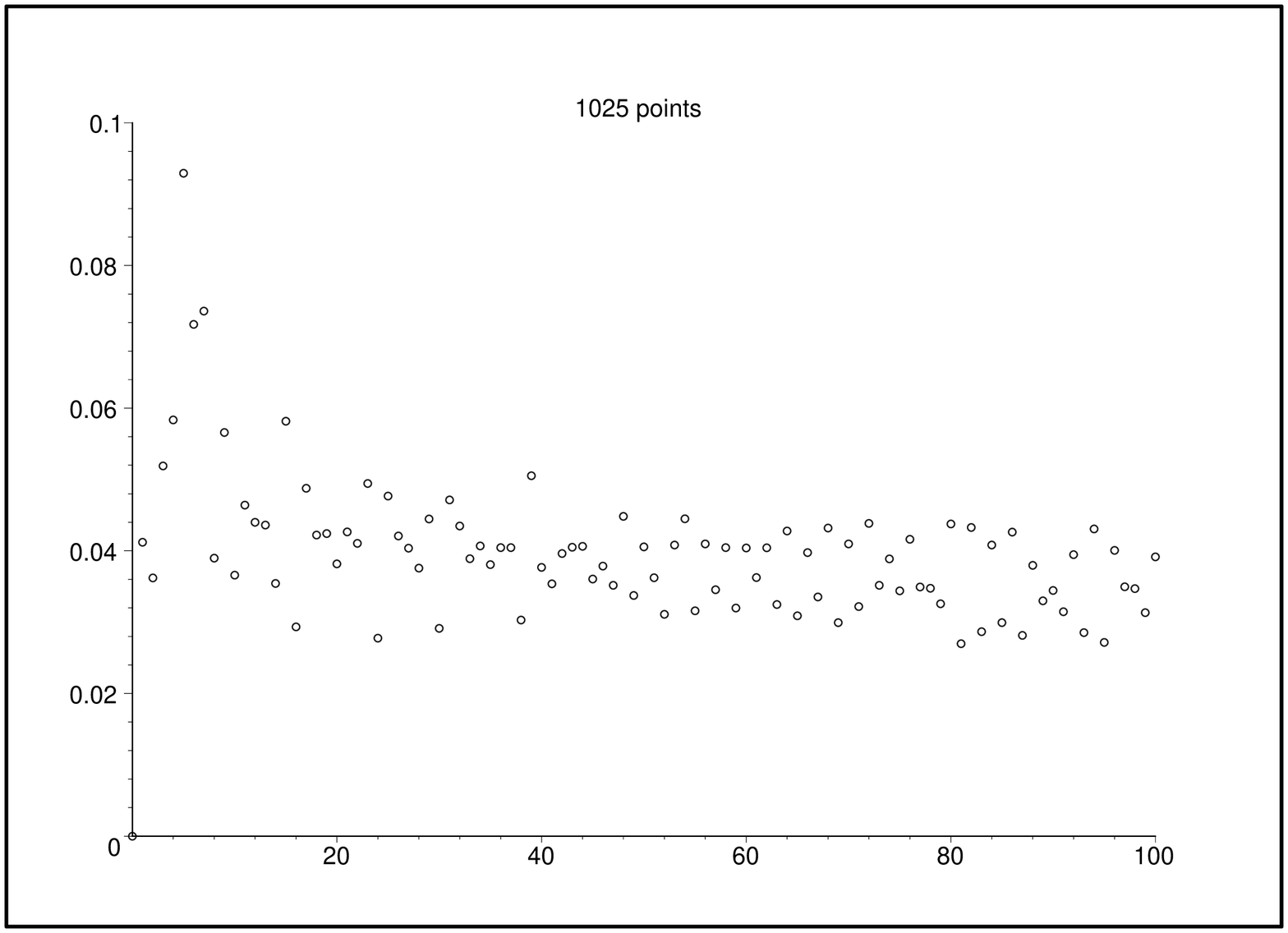}}}
\rotatebox{270}{\resizebox{2.3in}{!}{\includegraphics{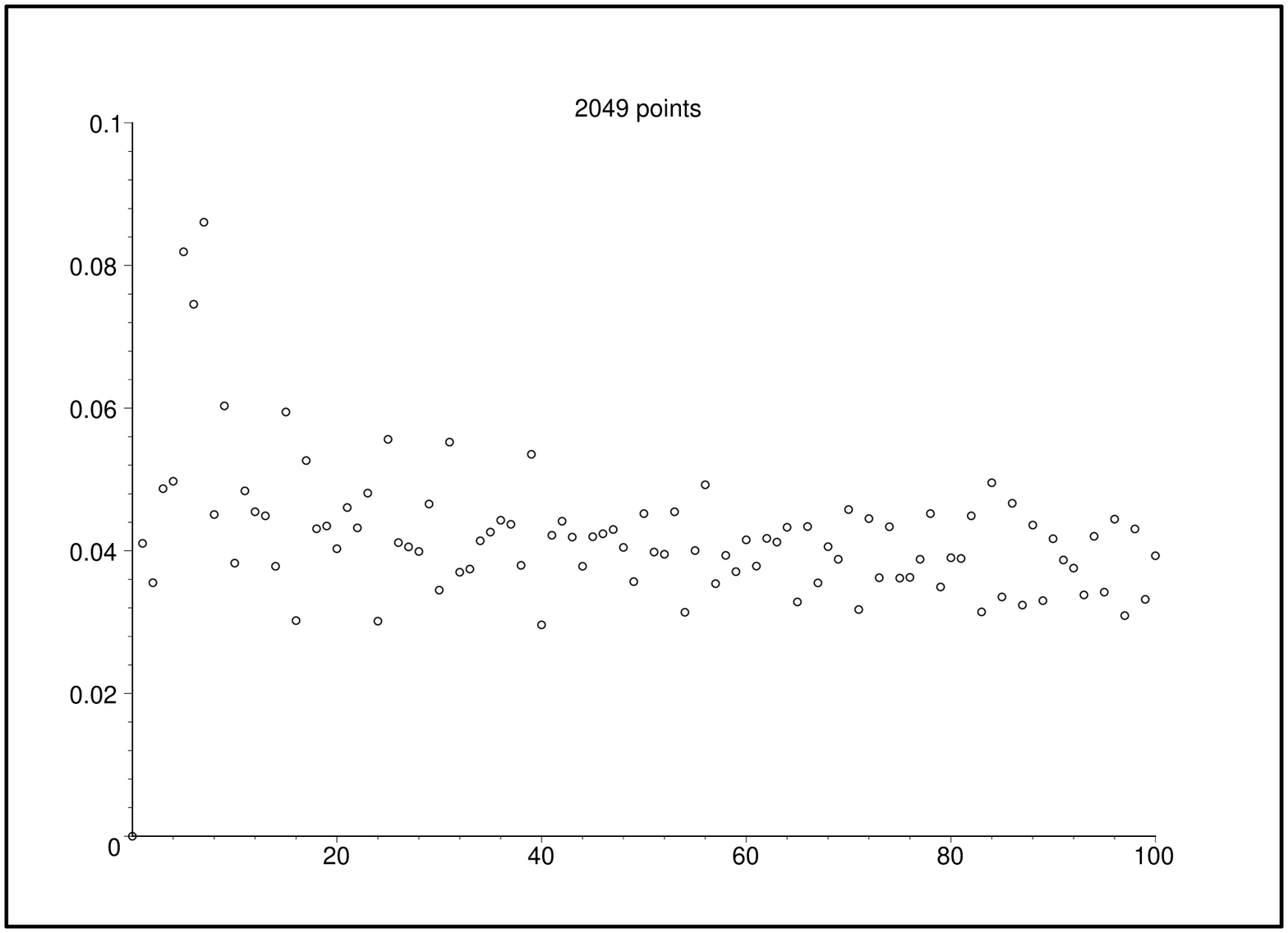}}}
\caption{The recurrence coefficients $a_{n,N}^2$ for $N=10, 11, 12$}
\label{fig101112}
\end{figure}

\begin{figure}[h!t]
\centering
\rotatebox{270}{\resizebox{2.3in}{!}{\includegraphics{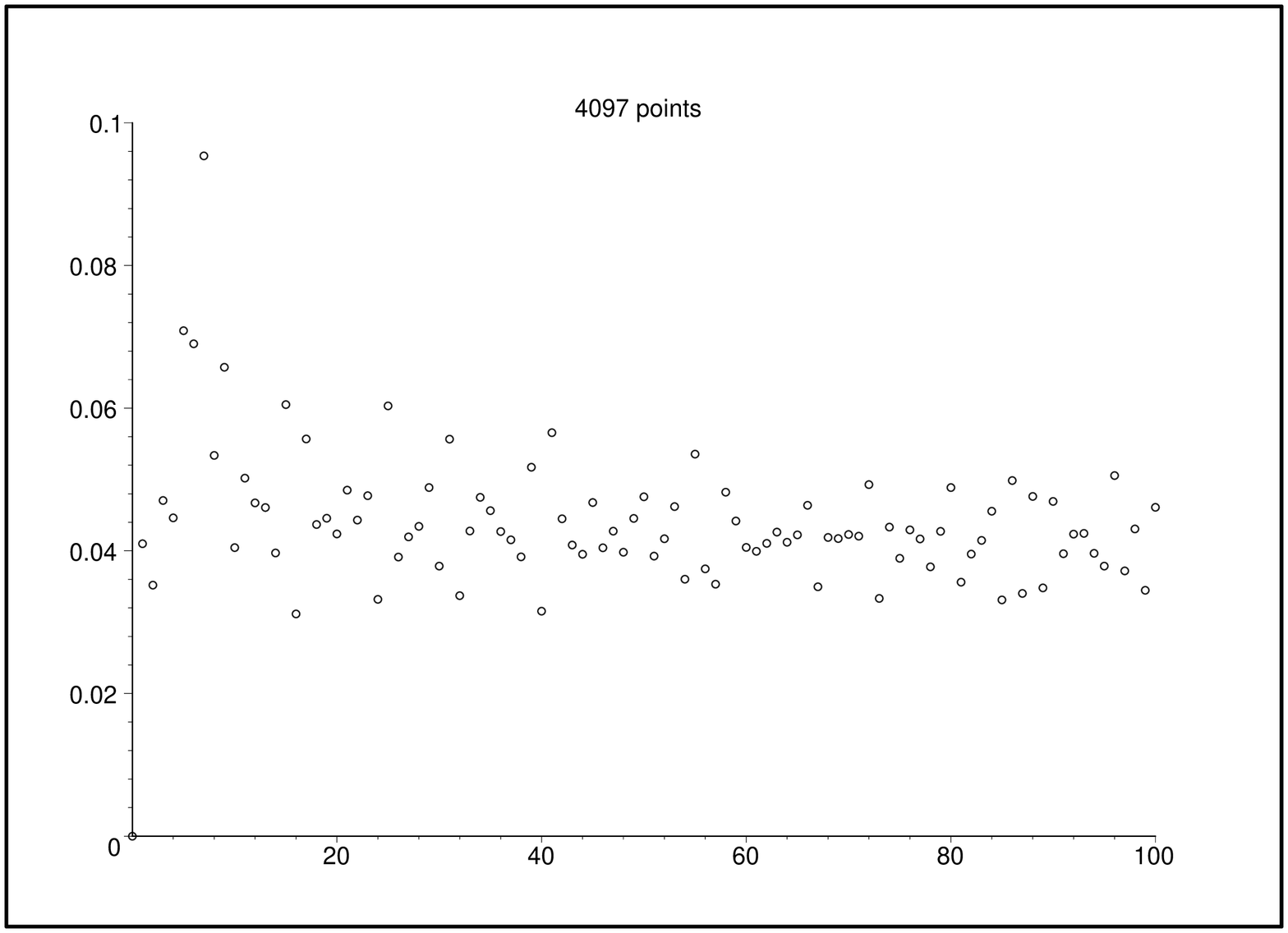}}}
\rotatebox{270}{\resizebox{2.3in}{!}{\includegraphics{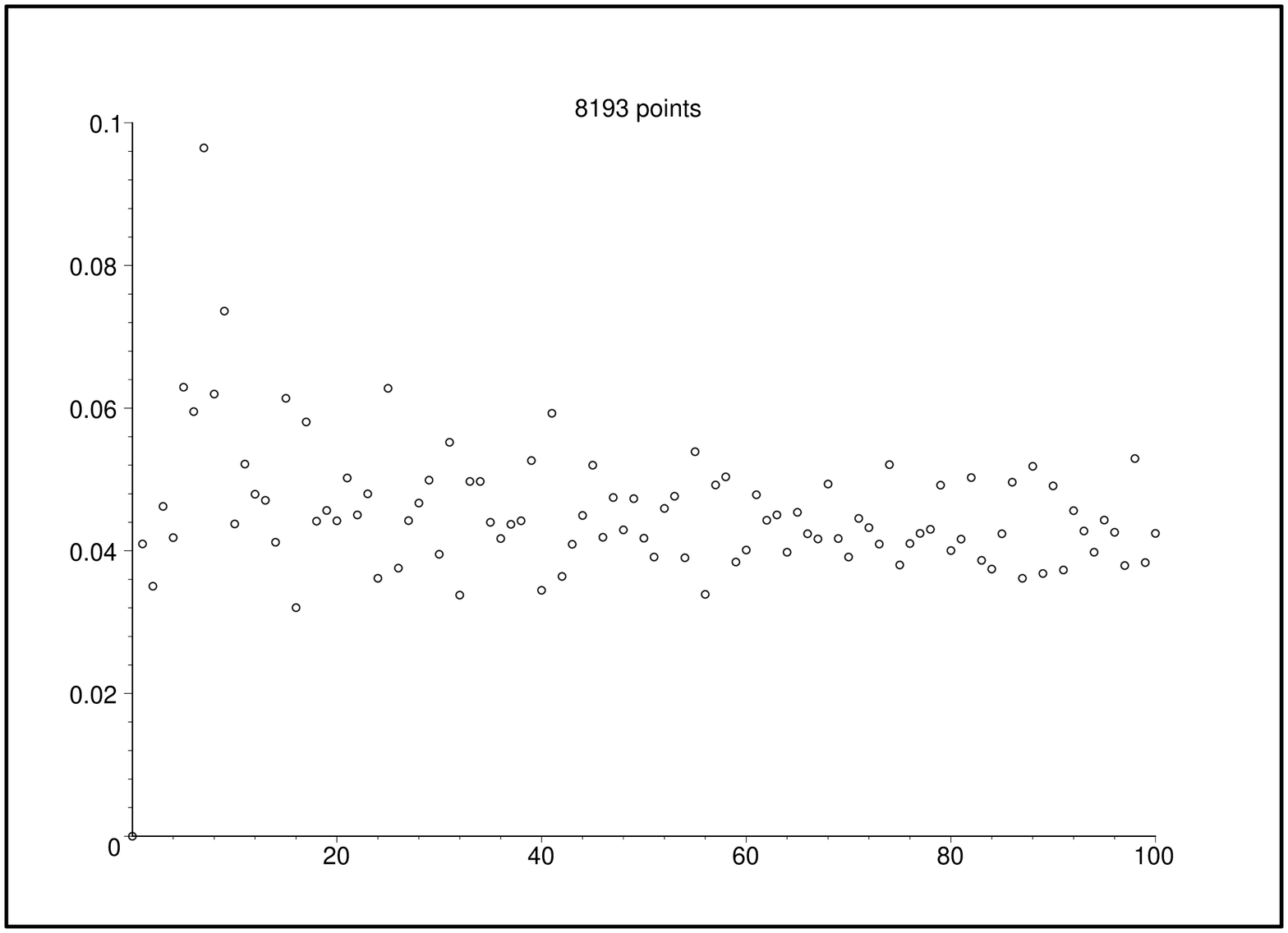}}}
\rotatebox{270}{\resizebox{2.3in}{!}{\includegraphics{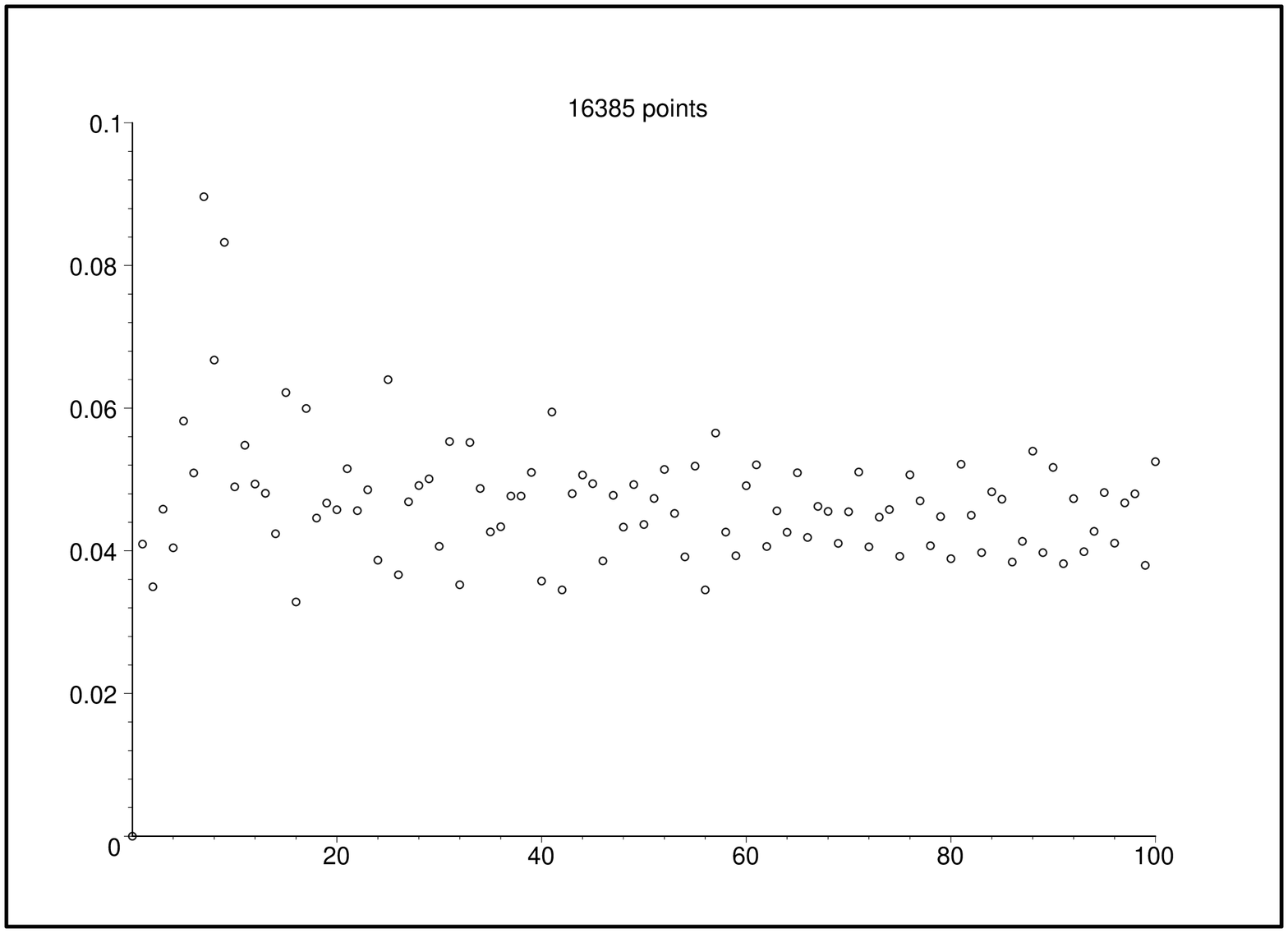}}}
\caption{The recurrence coefficients $a_{n,N}^2$ for $N=13, 14, 15$}
\label{fig131415}
\end{figure}

\begin{figure}[h!t]
\centering
\rotatebox{270}{\resizebox{2.3in}{!}{\includegraphics{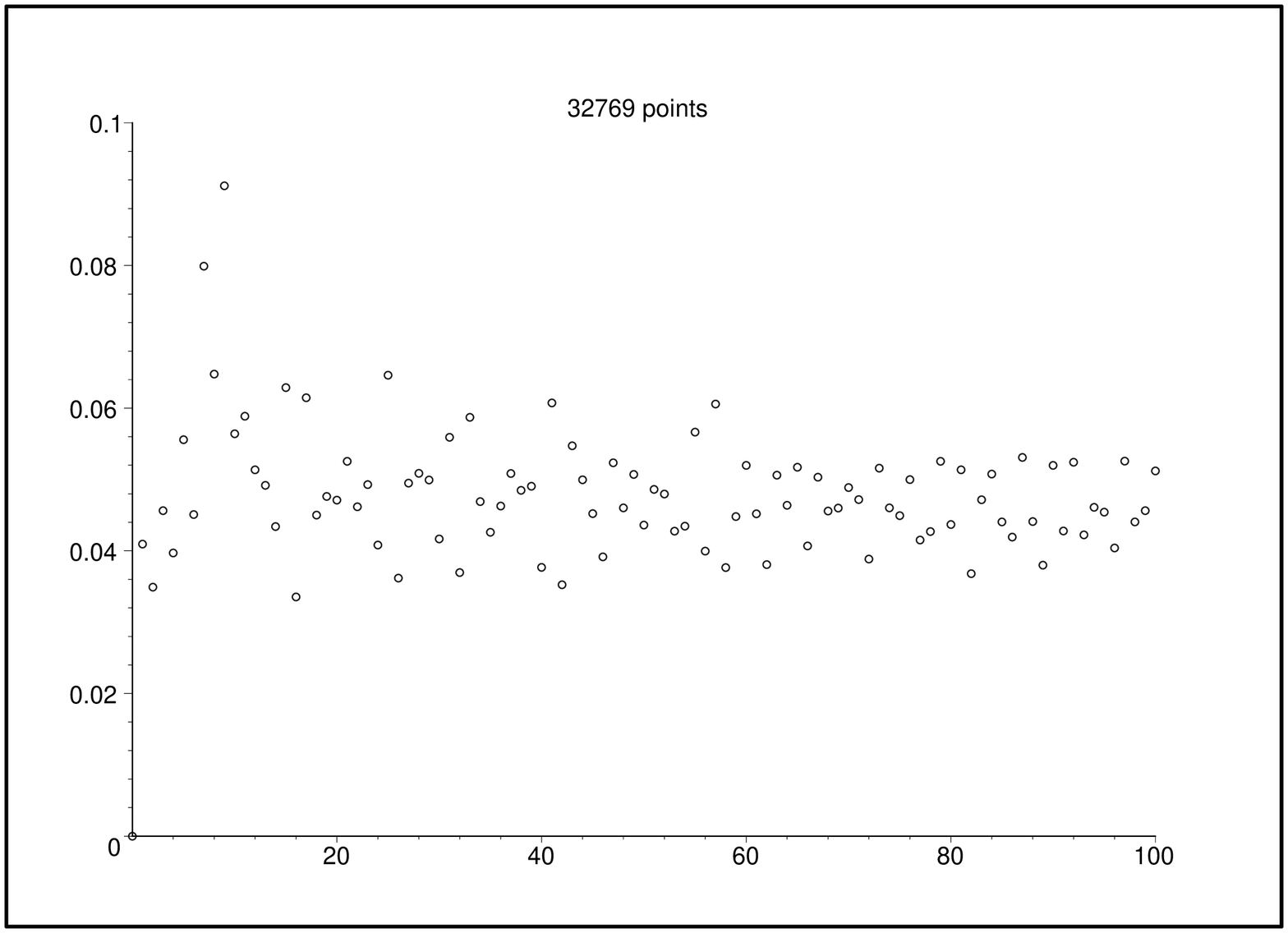}}}
\rotatebox{270}{\resizebox{2.3in}{!}{\includegraphics{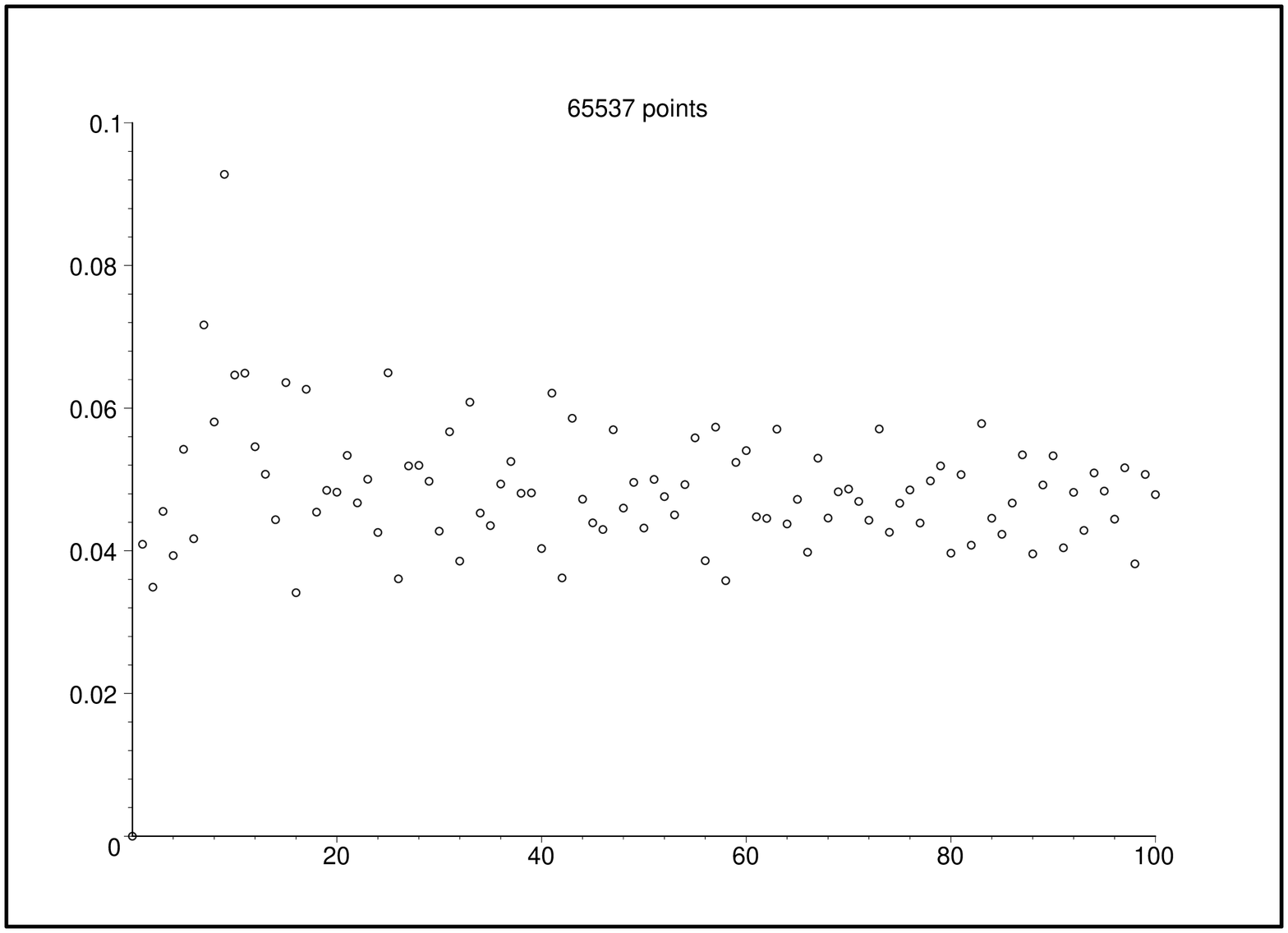}}}
\rotatebox{270}{\resizebox{2.3in}{!}{\includegraphics{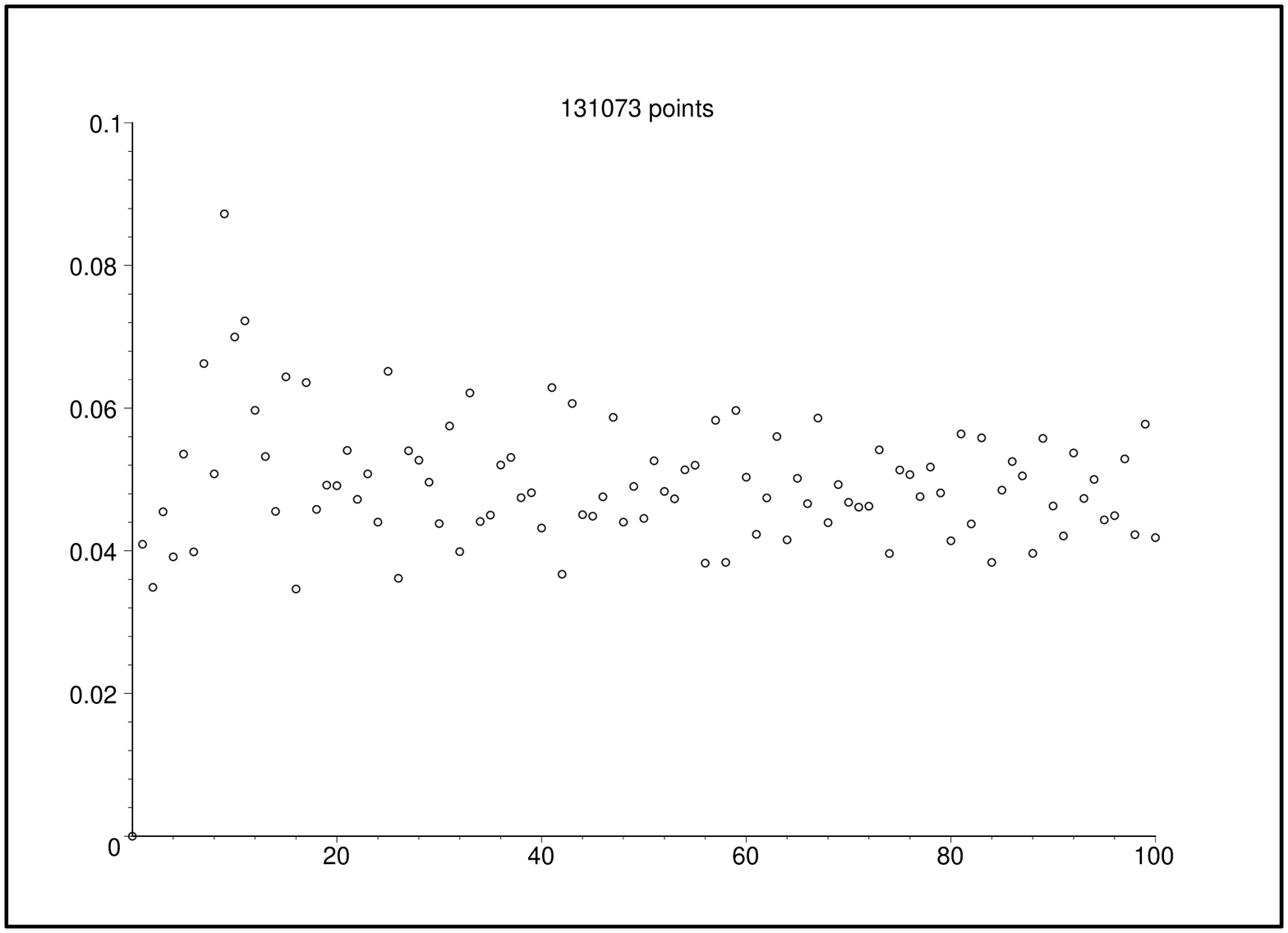}}}
\caption{The recurrence coefficients $a_{n,N}^2$ for $N=16, 17, 18$}
\label{fig161718}
\end{figure}

We have chosen to take for $q_N$ the empirical distribution function for the points in the Minkowski sequence $\mathcal{M}_N$. This means that we take
\[      q_N(x) = \frac{1}{m_N} \# \{ a/b \in \mathcal{M}_N \ | \ a/b \leq x \}, \]
where $m_N = 2^{N-1}+1$ is the number of points in $\mathcal{M}_N$. The measure induced by $q_N$ is supported on the Minkowski sequence
$\mathcal{M}_N$ and each point has equal mass $1/m_N$. The corresponding measure can be written as
\[   \frac{1}{m_N} \sum_{a/b \in \mathcal{M}_N} \delta_{a/b},  \]
where $\delta_c$ is the Dirac measure at $c$. It is not so difficult to see that this discrete measure converges weakly to the measure
induced by the Minkowski question mark function and in fact
\[   \lim_{N \to \infty}  q_N(x) = q(x), \]
uniformly on $[0,1]$ (see, e.g., \cite[\S 2.1]{Alkauskas}).
We have used these empirical distribution functions for $N=10$ up to $N=18$ (hence discrete measures with 513, 1025, 2049, 4097, 8193, 16385, 32769, 65537, and 131073 points) to compute the first 100 recurrence coefficients. The discrete measures are symmetric with respect to the point $1/2$, hence automatically $b_{n,N} = 1/2$
for $0 \leq n \leq 2^{N-1}$, so the main problem is to compute the recurrence coefficients $a_{n,N}^2$ and since we are only interested
in the limit for $N \to \infty$, we restrict our attention to $1 \leq n \leq 100$.
First we need to find all the points in the set $\mathcal{M}_N$, which is easily done by the following procedure in Maple, using the packages 
\texttt{LinearAlgebra} and \texttt{ArrayTools}:
\begin{verbatim}
minkseq:=proc(N)
  local s,k,i;
  s:=Vector[row](2^{N-1}+1,fill=1);
  s[1]:=0;
  s[2]:=1;
  for i from 2 to n do
     for k to 2^(i-2) do
        s[2^(i-2)+k+1]:=(numer(s[k])+numer(s[k+1]))
                        /(denom(s[k])+denom(s[k+1]))
     end do;
     s=:sort(s)
  end do;  
  s
end proc;
\end{verbatim}
The weights corresponding to the points in $\mathcal{M}_N$ are all equal to $1/(2^{N-1}+1)$. Then we used the discretized
Stieltjes-Gautschi method as is described in \cite[\S 4.1 on p.~34]{Gautschi1994} and implemented in the algorithm \texttt{stieltjes.m} \cite{OPQ}.
We have made our computations in Maple in a precision given by \texttt{Digits:=100}.
The resulting recurrence coefficients $a_{n,N}^2$ are given in Figures \ref{fig101112}--\ref{fig161718}.
\newpage
 
\newpage
 
\newpage

\newpage
 
Even though the figures all show a similar behavior, they give different approximations for the recurrence
coefficients of the orthogonal polynomials for the question mark function. The values are lying around $0.04$ to $0.05$ but there is a large value of $a_{n,N}^2$ (larger than 0.08)
near the beginning which changes significantly in size and position when $N$ changes, see Table \ref{large}.

\begin{table}
\centering
\begin{tabular}{|r|l|}
\hline
 $a_{5,10}^2$ & 0.0973813815 \\
 $a_{5,11}^2$ & 0.0929448256 \\ 
 $a_{7,12}^2$ & 0.0860780825 \\
 $a_{7,13}^2$ & 0.0953811702 \\
 $a_{7,14}^2$ & 0.0964882611 \\
 $a_{7,15}^2$ & 0.0896426756 \\
 $a_{9,16}^2$ & 0.0911621120 \\
 $a_{9,17}^2$ & 0.0927763299 \\
 $a_{9,18}^2$ & 0.0872304015 \\  
 \hline
\end{tabular}
\caption{Largest values of $a_{n,N}^2$}
\label{large}
\end{table}

\section{Moments}  \label{sec:moments}
In this section we will compute the recurrence coefficients by using the moments of the Minkowski question mark function.
We have computed the moments using a method suggested by G. Alkauskas \cite{Alkauskas1} \cite{Alkauskas2}. We use the relation
\begin{equation}  \label{mom}
   m_s = \sum_{k=0}^\infty (-1)^k c_{k+s} \binom{k+s-1}{k} m_k, \qquad s \geq 1, 
\end{equation}
where
\[   c_k = \sum_{n=1}^\infty \frac{1}{2^n n^k},  \]
then truncate the sum for $c_n$ to 400 terms and truncate the sum in \eqref{mom} to 500 terms, so that it becomes
a linear system of equations for the moments, with initial moment $m_0=1$. The error we make by truncating the sum for $c_k$ is 
\[     \sum_{n=401}^\infty \frac{1}{2^n n^k} < \sum_{n=401}^\infty \frac{1}{2^n} = 2^{-400} \approx 0.387 \ 10^{-120} , \]
hence the $c_k$ are accurately computed up to 120 decimals. 
In Maple we used \texttt{Digits:=400}
and we obtained the moments as given (with 50 decimals) in \ref{App:moments}.
The condition number of the matrix for this linear system is quite high (of the order $10^{435}$) so that the results must be
treated with some suspicion, even with using high accuracy.  
The first moment $m_1$ is 1/2 and hence this value can be used to check the accuracy of the method: the computed
value of $m_1$ was correct up to 93 decimals ($m_1=0.4999999\ldots$, with 92 nines). The accuracy seems to decrease for the
higher moments. Our value of $m_{100}$ is
\[   m_{100} = .000000444593386091498 \]
which differs slightly from the value $0.0000004445933003$ in \cite[p.~366]{Alkauskas1}.

An exact formula for the moments is given in \cite{Alkauskas3}. This formula is not very suitable for computing the moments, 
but another infinite system of equations was mentioned in \cite[p.~364]{Alkauskas4}. 
The relation is\footnote{We thank one of the referees for pointing this out}
\begin{equation}  \label{mom2}
   m_s = \sum_{k=0}^\infty d_{k+s} \binom{k+s-1}{k} m_k, \quad s \geq 1, 
\end{equation}
where
\[  d_k = 2 \sum_{n=2}^\infty \frac{1}{2^n n^k} = 2c_k-1.  \]
This relation can be proved in a similar way as Proposition 5 in \cite{Alkauskas3} but using another functional equation
for the moment generating function
\[  G(z) = \sum_{k=1}^\infty m_k z^{k-1}, \]
than the one used in the proof of \eqref{mom}, namely \cite[Eq.~(12)]{Alkauskas3}
\[  G(z) = - \frac{1}{1-z} - \frac{1}{(1-z)^2} G\left(\frac{1}{1-z}\right) + 2G(z+1).  \]
If we truncate the sum in \eqref{mom2} to 500 terms, then the linear system has a matrix with only positive terms and the
condition number 2.97 is low. The value of $m_1$ was correct up to $31$ decimals ($m_1=0.4999999\ldots$, with 30 nines), which is
less than using the linear system \eqref{mom}, but now this error is of the same magnitude for all moments. In fact we do get the same value
for the moment $m_{100}$ as given higher. We believe that the reason why
\eqref{mom2} gives less accurate results than \eqref{mom} is that the truncation of the infinite system of equations to a matrix of size 500
leaves the error
\[  \sum_{k=501}^\infty (-1)^k c_{k+s} \binom{k+s-1}{k} m_k \]
for the first system which, due to the oscillating terms in the sum, is much less than the error
\[  \sum_{k=501}^\infty d_{k+s} \binom{k+s-1}{k} m_k  \]
for the second system, which contains only positive terms. So even though the matrix of the second system is better conditioned than the matrix
of the first system, the results of the first system give more accurate approximations to the moments.

\begin{figure}[ht]
\centering
\rotatebox{270}{\resizebox{!}{5in}{\includegraphics{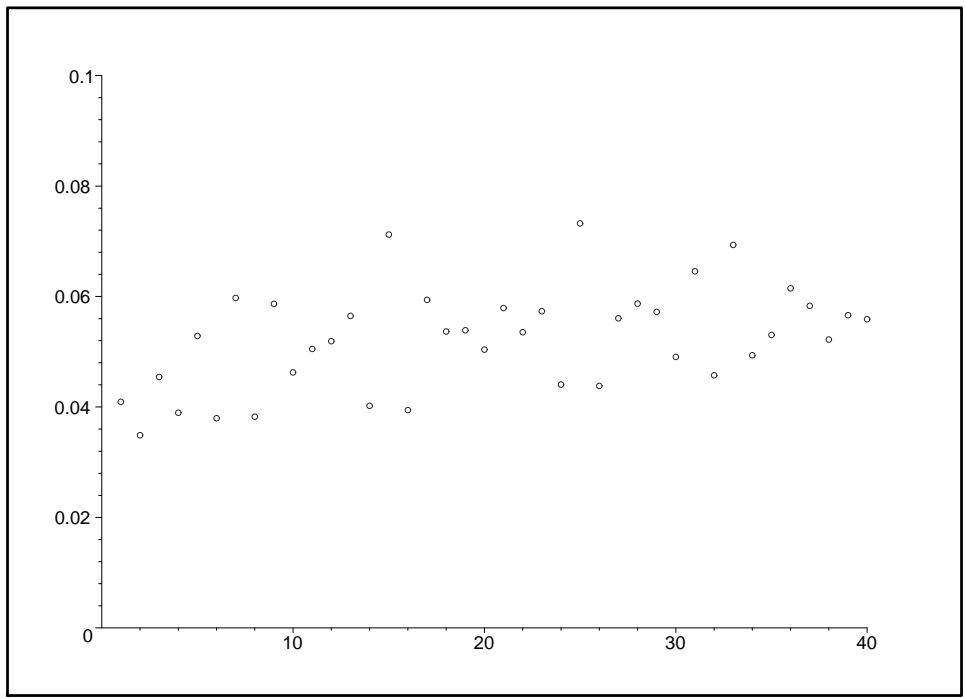}}}
\caption{The recurrence coefficients $a_n^2$}
\label{fig:an2}
\end{figure}

Next, we used these moments to compute the recurrence coefficients in the three-term recurrence relation
\[    P_{n+1}(x) = (x-b_n) P_n(x) - a_n^2 P_{n-1}(x)  \]
of the monic orthogonal polynomials for Minkowski's question mark function. We used the Chebyshev algorithm (with the ordinary moments, see \cite[Algorithm 2.1 on p.~77]{Gautschi}). It is well known that the mapping from moments $(m_n)_{n \geq 0}$ to the recurrence coefficients
$(a_{n+1}^2,b_n)_{n\geq 0}$ is badly conditioned, see, e.g., \cite[\S 2.1.6]{Gautschi}. This is the reason why we used high precision 
(\texttt{Digits:=400}) in our calculations. 
Recall that, due to the symmetry, all the recurrence coefficients $b_n$ are constant: $b_n = 1/2$. This was useful to check the accuracy of our computations. We observed that our computed values of $b_n$ were correct to 23 decimal places up to $b_{40}$ but then slowly started to show errors, 
with $b_{53}$  a negative value, which is impossible. Therefore we listed and plotted the computed values of $a_n^2$ only up to
$a_{40}^2$ in Table \ref{table-ab} and Figure \ref{fig:an2}. 

\begin{table}
\vspace{-1cm}
\centering
\small
\begin{tabular}{|r|l|l|}
\hline
 $k$ & $b_k$  & $a_k^2$ \\
 \hline
 0 & 0.5 &  \\
 1 & 0.5 & 0.040926476429308736381 \\
 2 & 0.5 & 0.034881265506134342903 \\
 3 & 0.5 & 0.045430415370805808038 \\
 4 & 0.5 & 0.038973377115288248098 \\ 
 5 & 0.5 & 0.052863907245188596784 \\
 6 & 0.5 & 0.037955175327144719607 \\
 7 & 0.5 & 0.059731637094523918352 \\
 8 & 0.5 & 0.038238400877758730555 \\ 
 9 & 0.5 & 0.058672115522904960765 \\ 
 10 & 0.5 & 0.046255346737208862213 \\
 11 & 0.5 & 0.050520824434494850803 \\ 
 12 & 0.5 & 0.051910925145095030363 \\ 
 13 & 0.5 & 0.056489563093038456301 \\ 
 14 & 0.5 & 0.040208992500495472293 \\ 
 15 & 0.5 & 0.071218137450141992615 \\ 
 16 & 0.5 & 0.039427602611174900647 \\ 
 17 & 0.5 & 0.059396186789821055700 \\ 
 18 & 0.5 & 0.053652031489601189600 \\ 
 19 & 0.5 & 0.053884790282064402379 \\ 
 20 & 0.5 & 0.050381653151077022836 \\ 
 21 & 0.5 & 0.057911359198380156348 \\ 
 22 & 0.5 & 0.053527412587600219313 \\
 23 & 0.5 & 0.057334758849746482997 \\ 
 24 & 0.5 & 0.044067432839352949172 \\
 25 & 0.5 & 0.073234222409597016726 \\ 
 26 & 0.5 & 0.043818059541906812748 \\  
 27 & 0.5 & 0.056063579773800371687 \\  
 28 & 0.5 & 0.058703247843561897668 \\  
 29 & 0.5 & 0.057201243688393241195 \\  
 30 & 0.5 & 0.049069569275476665894 \\  
 31 & 0.5 & 0.064554576616699275413 \\
 32 & 0.5 & 0.045732073443752859070 \\
 33 & 0.5 & 0.069343850504209381053 \\
 34 & 0.5 & 0.049351815705867268360 \\
 35 & 0.5 & 0.053031677464673738708 \\
 36 & 0.5 & 0.061496568432752989923 \\
 37 & 0.5 & 0.058321435643186334303 \\
 38 & 0.5 & 0.052208175547033955364 \\
 39 & 0.5 & 0.056607898016942758422 \\
 40 & 0.5 & 0.055895931809777873999 \\
\hline
\end{tabular}
\caption{The recurrence coefficients $b_n$ and $a_n^2$}
\label{table-ab}
\end{table}

\newpage

In order to check whether the question mark function induces a regular measure on $[0,1]$, we need to compute $\gamma_n^{-1/n}$, where 
\[  \frac{1}{\gamma_n^2} = \int_0^1 P_n^2(x)\, dq(x) = a_1^2 a_2^2 \cdots a_n^2 . \]
For a regular measure on $[0,1]$ one has
\[   \lim_{n \to \infty} \left( \frac{1} {\gamma_n} \right)^{1/n} = \frac14, \]
where $1/4$ is the logarithmic capacity of the interval $[0,1]$. It is more convenient to compute $\gamma_n^{-2/n}$ since we are really computing the squared recurrence coefficients $a_n^2$. Our computations (see Table \ref{regular?}) indicate that the sequence
$( a_1^2a_2^2\cdots a_n^2)^{1/n}$ may indeed be converging, but the limit is near $0.052$ instead of $(1/4)^2 = 1/16 = 0.0625$.

\begin{table}
\vspace{-1cm}
\centering
\small
\begin{tabular}{|r|l|}
\hline
 $k$ & $(a_1^2a_2^2\cdots a_k^2)^{1/k}$ \\
\hline
 1 & 0.040926476429308736381 \\
 2 & 0.037783161468586334476 \\
 3 & 0.040177332455719534022 \\
 4 & 0.039872900960755211548 \\
 5 & 0.042186560554634906505 \\
 6 & 0.041449910727599868208 \\
 7 & 0.043670914641443765852 \\ 
 8 & 0.042951735555887041699 \\
 9 & 0.044466283404646215085 \\
 10 & 0.044642030777323595493 \\
 11 & 0.045146924287480657337 \\
 12 & 0.045675227441320883038 \\
 13 & 0.046427976597463213256 \\ 
 14 & 0.045953497782193798426 \\
 15 & 0.047315493717202662999 \\
 16 & 0.046779242858071064616 \\
 17 & 0.047440963915026139792 \\
 18 & 0.047766342240858080778 \\
 19 & 0.048070312391450962315 \\
 20 & 0.048183319646452786085 \\
 21 & 0.048607122313805874900 \\ 
 22 & 0.048820630103435556238 \\ 
 23 & 0.049163047752620715919 \\
 24 & 0.048939412842737333783 \\ 
 25 & 0.049734867738244766161 \\
 26 & 0.049493171229376460057 \\  
 27 & 0.049722196075856275919 \\  
 28 & 0.050017931134405002549 \\  
 29 & 0.050249919542639211311 \\  
 30 & 0.050210120837211061076 \\  
 31 & 0.050618791840573811015 \\
 32 & 0.050458453441009796727 \\
 33 & 0.050946927334564490546 \\
 34 & 0.050899284439881664328 \\
 35 & 0.050959003303570963593 \\
 36 & 0.051225761593643949404 \\
 37 & 0.051405681461085926182 \\
 38 & 0.051426640855254603168 \\
 39 & 0.051553375093489825709 \\
 40 & 0.051657713625567224389 \\
 \hline
\end{tabular}
\caption{The values of $\gamma_n^{-2/n}=(a_1^2a_2^2\cdots a_n^2)^{1/n}$.}
\label{regular?}
\end{table}

\section{Conclusion}
\subsection{The discretized Stieltjes-Gautschi method}
When we compare the values of $a_{n,N}^2$ obtained in Section 2 with the values obtained in Section 3 (using the moments of $q$), 
see  Figure \ref{fig:an2}, only the first five recurrence coefficients 
$a_{n,18}^2$ are close to the actual values. The reason that we don't get accurate results for the $a_n^2$ is that the $a_{n,N}^2$ are
approximations for the $a_n^2$ and that one needs to take the limit $N \to \infty$ to get the desired recurrence coefficients.
The discrete distribution functions $q_N$ converge slowly to the Minkowski question mark function. In fact, it is known that
\[   \max_{x \in [0,1]} |q_N(x)-q(x)| = \frac{1}{2^{N-1}+1}, \]
hence the error $\|q_N-q\|_\infty$ is of the order $1/m_N$, where $m_N$ is the number of points in the support of the discrete measure $q_N$.
The discrete approximation $q_N(x)$ for $x$ near $0$ and near $1$ is quite poor since the first point in $\mathcal{M}_N$ (after $0$) is $1/N$
and the largest point (before 1) is $(N-1)/N$. These single points have the task to represent $q$ on an interval of reasonable large size compared to the number of points which are available. The computations are in fact quite accurate for the recurrence coefficients $a_{n,N}^2$ of the discrete orthogonal polynomials, but the convergence in \eqref{discreteapprox} is rather slow (except, of course, for $b_{n,N}$).    

\subsection{Behavior of the recurrence coefficients}
If we compute the recurrence coefficients $(a_{k}^2)_{1 \leq k \leq n}$ and $(b_k)_{0 \leq k \leq n}$ from the moments $(m_k)_{0\leq k \leq 2n}$, then we can't compute many coefficients since the mapping from moments to recurrence coefficients becomes ill-conditioned at an exponential rate in $n$. 
The use of modified moments 
\[   \hat{m}_n = \int_0^1 r_n(x)\, dq(x),  \]
where $(r_n)_{n \in \mathbb{N}}$ is a sequence of known polynomials, sometimes leads to a better conditioned mapping from modified moments
to recurrence coefficients, in particular when the polynomials $(r_n)_{n \in \mathbb{N}}$ are already close to the polynomials $(p_n)_{n \in \mathbb{N}}$
for which we are computing the recurrence coefficients. Unfortunately no such system of polynomials seems to be available
and our attempts to use the Chebyshev polynomials on $[0,1]$ lead to a similar ill-conditioned problem. So we are stuck with the regular moments
and high precision arithmetic, which allowed us to compute 40 recurrence coefficients $a_n^2$ reliably. The reliability was checked by using
the recurrence coefficients $b_k$ as control variables, since we know that they are all constant and equal to 1/2: errors in the computed values of $b_k$
for $k > n_0$ certainly indicate that the corresponding $a_k^2$ for $k > n_0$ are not reliable. On the other hand, if the $b_k$ are computed accurately
for $k \leq n_0$ then one may reasonably assume that the $a_k^2$ are also accurate for $k \leq n_0$, except possibly the last one. 

Of course, 40 recurrence coefficients are not really enough to say something about the asymptotic behavior of the recurrence coefficients.
Nevertheless Figure \ref{fig:an2} doesn't give the impression that the recurrence coefficients are converging. The recurrence coefficients vary somewhat
between $0.04$ and $0.06$ with some higher values. The average of the recurrence coefficients $a_n^2$ for $1 \leq 40$ is  $0.05246234283$ which is not
close to $1/16=0.0625$. Hence, based on our numerical evidence, we conclude that the recurrence coefficients $(a_n^2)_n$ do not converge, and even their averages do not seem to converge to $1/16$. This means that Minkowski's question mark function \textbf{does not give orthogonal polynomials in the Nevai class}
$M(\frac12,\frac14)$, and most likely not in any other Nevai class $M(\frac12,a)$.    

\subsection{Does $q$ induce a regular measure?}
The oscillating nature of the recurrence coefficients $(a_n^2)_{n \geq 1}$ still leaves open the possibility that the geometric mean
$(a_1^2a_2^2\cdots a_n^2)^{1/n}$ converges. Our numerical experiments, with $1\leq n \leq 40$, indicate that the geometric mean may be converging
but that the limit is less than $1/16=0.0625$, see Table \ref{regular?}. Hence our numerical evidence leads to the conclusion that the question mark function  \textbf{does not induce a regular measure} on $[0,1]$.

In general it is known that
\[     c_\mu   \leq \liminf_{n \to \infty}  \gamma_n^{-1/n} \leq  \limsup_{n\to \infty} \gamma_n^{-1/n}  \leq  \textup{cap}(S_\mu), \]
where $S_\mu$ is the support of the orthogonality measure $\mu$ and $c_\mu$ is the minimal carrier capacity of $\mu$,
\[   c_\mu = \inf \{ \textup{cap}(B)\ | \ B \textrm{ is a bounded Borel set and } \mu(\mathbb{R} \setminus B) = 0\}.  \]
Hence if the Minkowski question mark function does not induce a regular measure $\mu$ on $[0,1]$, then this implies that the minimal carrier
capacity of the question mark function is less than $1/4$. The reason why $q$ does not induce a regular measure is probably that the support
$[0,1]$ contains intervals with exponentially small measure. Indeed, one easily finds
\[     q\left(\frac1{n}\right)-q(0) = \frac{1}{2^{n-1}}, \qquad q(1) - q\left(1-\frac{1}{n}\right) = \frac{1}{2^{n-1}}, \]
\[     q\left(\frac12+\frac{1}{4n}\right) - q\left(\frac12-\frac{1}{4n}\right) = \frac32 \frac{1}{2^n}, \]
so that the intervals $[0,\frac{1}{n}]$, $[\frac{n-1}{n},1]$ and $[\frac12-\frac1{4n},\frac12+\frac1{4n}]$ of length $\mathcal{O}(1/n)$ have measure proportional to $1/2^n$. These intervals reappear all over the interval $[0,1]$
because of the self-similar nature of $q$, expressed by \eqref{self-similar}.
Hence, the support of the measure induced by the question mark function behaves like a subset of $[0,1]$ with gaps and for measures
with a support containing gaps, the recurrence coefficients show much more oscillating or chaotic behavior. The minimal carrier capacity
would be less than the capacity of $[0,1]$ due to these gaps which have exponentially small measure.

\subsection{Open problems}
Many open problems remain for the orthogonal polynomials related to the question mark function \cite{open}. To name a few:
\begin{itemize}
\item Is there any structure in the seemingly chaotic behavior of the recurrence coefficients? In particular it would be nice to know whether
the coefficients $a_n^2$ are almost periodic. To answer this question it would be necessary to calculate the mean of the coefficients and then to look for oscillations about this mean. This would require many more coefficients than available right now. 
\item Can one determine
\[   \liminf_{n \to \infty} a_n^2, \quad  \limsup_{n \to \infty}  a_n^2 \ ? \]
\item Does the geometric mean 
\[   \lim_{n \to \infty}  (a_1a_2\ldots a_n)^{2/n}   \]
exist? 
\item Is this limit (if it exists) equal to the minimal carrier capacity of the measure induced by the question mark function?  
\end{itemize}

The question mark function is a reasonably simple singularly continuous measure which, next to the well known Cantor measure, is useful to investigate
how the recurrence coefficients of orthogonal polynomials with singular continuous measure behave. This paper is the first in which the
Minkowski question mark function is considered in the context of orthogonal polynomials. Earlier there have been a number of papers
dealing with the recurrence coefficients of orthogonal polynomials for the Cantor measure, and for these many more recurrence coefficients
have been investigated, see e.g., the work of G. Mantica \cite{Mantica}, H.-J. Fischer \cite{Fischer} and the recent work of 
Heilman, Owrutsky and Strichartz \cite{HOS} on orthogonal polynomials with self-similar measures. 

A very relevant singular continuous measure is the equilibrium measure for the Julia set of the iteration of a polynomial $T$, such as $T(x)=x^2-c$, with $c>2$. Such a Julia set is of the same nature as the Cantor set, where one iteratively removes intervals from a given interval. The orthogonal polynomials for such a singular measure have been analyzed in detail by Barnsley, Geronimo, Harrington \cite{BGH}, Bellissard, Bessis, Moussa \cite{BBM}, Bessis, Geronimo, Moussa \cite{BGM} and Bessis, Mehta, Moussa \cite{BMM}. One of their results is that the subsequence $P_{2^n}(x)$ of the orthogonal polynomials is explicitly given by the 
$n$-th iterate of the given polynomial $T(x)$ and that the recurrence coefficients satisfy some non-linear relations, from which one can deduce
that the recurrence coefficients are limit periodic for a large class of polynomials $T$. It would be nice to obtain such results for the Minkowski question mark function. However, the self-similarity of the question mark function, as described by \eqref{self-similar}, involves rational functions
so that (orthogonal) polynomials composed with rational functions lead to rational functions and hence the polynomial nature is not preserved by the mappings in the iterated function system.

\subsection*{Acknowledgement}
We thank S. Yakubovich for bringing the Minkowski question mark function to our attention.
This research was supported by KU Leuven research grant OT/12/073, FWO grant G.0427.09 and the Belgian Interuniversity Attraction Poles
Programme P7/18.

\newpage

\appendix

\section{Moments}   \label{App:moments}

\begin{tabular}{|r|l|}
\hline
 $k$ & $m_k$ \\
 \hline
 1 &  0.50000000000000000000000000000000000000000000000000 \\
 2 &  0.29092647642930873638069776273912029008043710219559 \\
 3 &  0.18638971464396310457104664410868043512065565329339 \\
 4 &  0.12699225840744313520289222788021163884118514576173 \\
 5 &  0.090164454945335997055486162852728371901870108915328 \\
 6 &  0.065928162577472341926815287122643139097250419686219 \\
 7 &  0.049294310463767495715873019738162192824596016540408 \\
 8 &  0.037518711852048335539887123885001326828330311845185 \\
 9 &  0.028979622034097514125202534817631105526453454918480 \\
 10 & 0.022665858176292038567084595433830909709028113615601 \\
 11 & 0.017920859234922559709620674061965959387422365122601 \\
 12 & 0.014304689510828028713327933529010923388486775628519 \\
 13 & 0.011515014023688037216614175868946822321198574437066 \\
 14 & 0.0093396445516456630408229433763950525836249376062925 \\
 15 & 0.0076269557590972324137284663579204452441313886605484 \\
 16 & 0.0062668792729955855181245666105463942184539549010276 \\
 17 & 0.0051783877867685200880408849955806484581000333322566 \\
 18 & 0.0043010785465847754710770555353762470229805125766226 \\
 19 & 0.0035894093787684071613366998806631220004089279356802 \\
 20 & 0.0030086867071149232599754021703205599242915788661537 \\
 21 & 0.0025322317634222266277036639307882264778052530156116 \\
 22 & 0.0021393519928752509332902769549514723550086345056961 \\
 23 & 0.0018138708773921687964413975361770115813330394837942 \\
 24 & 0.0015430503000888974750639610156922661098154179170794 \\
 25 & 0.0013167923220646751068501044957069765425398359511571 \\
 26 & 0.0011270421715775897395157816468683456633974505562164 \\
 27 & 0.00096733770677037783561122437105791592882537796051057 \\
 28 & 0.00083246658220639104473557479405732069610975267509581 \\
 29 & 0.00071820335559030841469949022408376381202557302200995 \\
 30 & 0.00062110644645096567556056888154034901542203410280702 \\
 31 & 0.00053836027103037905202629197075705905383821217884329 \\
 32 & 0.00046765173455275558510114895133989297915712195884076 \\
 33 & 0.00040707303776263157571881765213893343503728847982982 \\
 34 & 0.00035504477084132941930948314931026013893192071414219 \\
 35 & 0.00031025474516371162917662025303647309624750426315414 \\
 \hline
\end{tabular}

\begin{tabular}{|r|l|}
\hline
 $k$ & $m_k$ \\
 \hline
 36 & 0.00027160910474659558624707454363471603060022780445456 \\
 37 & 0.00023819307170515894210646197643373870213103110230062 \\
 38 & 0.00020923928922626913333968336609713458623759072879973 \\
 39 & 0.00018410218544426637502507882484141316037994511212385 \\
 40 & 0.00016223713098006110726991119269902592433404152806255 \\
 41 & 0.00014318342992835510310721778791488526305144139896956 \\
 42 & 0.00012655038932266879782766519282463067535339428877894 \\
 43 & 0.00011200587071717784852138503115054268313120729105696 \\
 44 & 0.000099266850721688056248146847039352100867092413446553 \\
 45 & 0.000088091613483408666548718142966352919966384794835161  \\  
 46 & 0.000078273273509573469871441756231586836413236191991855  \\  
 47 & 0.000069634386611004882448594820251627462598040246563183  \\  
 48 & 0.000062022453717890701051430866178792091112759697349794  \\  
 49 & 0.000055306159621487594284446897034200886840998808601394  \\  
 50 & 0.000049372218434717208976312359496974514537772757968773  \\  
 51 & 0.000044122721362726073881354844160254923500862547763273  \\  
 52 & 0.000039472901486704070959737735999458212044994283495032  \\  
 53 & 0.000035349245666258223798489554957594867788549703093561  \\  
 54 & 0.000031687896118933003424939851849986040400796168900580  \\  
 55 & 0.000028433294336821167964517994319875407187553962659664  \\  
 56 & 0.000025537028219094604542709782295425670986358747374410  \\  
 57 & 0.000022956850006417952618424640438139791507342597816299  \\  
 58 & 0.000020655838092375878340873282295400282104460092785749  \\  
 59 & 0.000018601680291880423863588358493197745228211014253879  \\  
 60 & 0.000016766059853447792758997531437488191617648704513657  \\  
 61 & 0.000015124128560475439274975373678569160508643322805383  \\  
 62 & 0.000013654053796046797396039913645973219220649444564075  \\  
 63 & 0.000012336628542866380312339531772945578083278145703849  \\  
 64 & 0.000011154935032657469162573155131448984539990112465750  \\  
 65 & 0.000010094054210908289421180310562947718378648386875946  \\  
 66 & 0.0000091408143945475299648985162271152064510786088610998 \\  
 67 & 0.0000082835735137657351740548705254404013142671182635503 \\  
 68 & 0.0000075120301789106204917387855457822400686084713220510 \\  
 69 & 0.0000068170595271178220756077387776286434567656799242417 \\  
 70 & 0.0000061905704040288945388332526110642901215894448069398 \\  
\hline
\end{tabular}

\end{document}